\documentclass[12pt,reqno]{amsart}
\usepackage{amsmath,amssymb,mathrsfs,amsthm,amsfonts}
\usepackage[usenames,dvipsnames]{xcolor}
\usepackage[inline]{enumitem} 
\usepackage{hyperref}
\hypersetup{%
	colorlinks=true, linkcolor=blue, 
	citecolor=ForestGreen
}

\usepackage[paper=letterpaper,margin=1in]{geometry}

\newtheorem{theorem}{Theorem}[section]

\newtheorem{corollary}[theorem]{Corollary}
\theoremstyle{definition}

\newtheorem{remark}[theorem]{Remark}

\numberwithin{equation}{section}
\allowdisplaybreaks

\usepackage{acronym}
\acrodef{SDE}{Stochastic Differential Equation}

\newcommand{\Ex}{ \mathbf{E} }				
\renewcommand{\Pr}{ \mathbf{P} }			
\newcommand{\ind}{\mathbf{1}}				
\newcommand{\set}[1]{{\{#1\}}}				

\newcommand{\Exp}{\operatorname{Exp}}
\newcommand{\Gammad}{\operatorname{Gamma}}
\newcommand{\Pois}{\operatorname{Pois}}
\newcommand{\supp}{\operatorname{supp}}

\newcommand{\per}{p}			

\newcommand{\poi}{\mu}			
\newcommand{\invt}{\nu}			

\newcommand{\e}{\varepsilon}	

\newcommand{\N}{\mathbb{N}}
\newcommand{\R}{\mathbb{R}}

\newcommand{\calU}{\mathcal{U}}
\newcommand{\calV}{\mathcal{V}}
\newcommand{\calW}{\mathcal{W}}

\newcommand{\phis}{\phi^\text{s}}

\newcommand{\Cc}{C_\text{c}}

\renewcommand{\bar}{\overline}

\usepackage{graphicx}
\newcommand*{\Cdot}{{\raisebox{-0.5ex}{\scalebox{1.8}{$\cdot$}}}} 

\begin{document}

\title[Stationary Distributions of the Atlas Model]
{Stationary Distributions of the Atlas Model}
\author[L.-C.\ Tsai]{Li-Cheng Tsai}
\address{L.-C.\ Tsai,
	Departments of Mathematics, Columbia University,
	\newline\hphantom{\quad\quad L-C Tsai}
	2990 Broadway, New York, NY 10027}
\email{lctsai.math@gmail.com}
\subjclass[2010]{
Primary 60J60, 	
Secondary 60H10.	
}

\maketitle

\begin{abstract}
In this article we study the Atlas model,
which consists of Brownian particles on $ \R $,
independent except that the Atlas (i.e., lowest ranked) particle $ X_{(1)}(t) $
receives drift $ \gamma dt $, $ \gamma\in\R $.
For any fixed shape parameter $ a>2\gamma_- $,
we show that, up to a shift $ \frac{a}{2}t $,
the \emph{entire} particle system has an invariant distribution $ \invt_a $,
written in terms an explicit Radon-Nikodym derivative
with respect to the Poisson point process of density $ ae^{a\xi} d\xi $.
We further show that $ \invt_a $ indeed has the product-of-exponential
gap distribution $ \pi_a $ derived in \cite{sarantsev17}.
As a simple application,
we establish a bound on the fluctuation of the Atlas particle $ X_{(1)}(t) $ uniformly in $ t $,
with the gaps initiated from $ \pi_a $ and $ X_{(1)}(0)=0 $.
\end{abstract}

\section{Introduction}
\label{sect:intro}

In this article we study the (infinite) Atlas model.
Such a model consists of
a semi-infinite collection of particles $ X_i(t) $, $ i=1,2,\ldots $,
performing independent Brownian motions on $ \R $,
except that the Atlas (i.e., lowest ranked) particle receives a drift of strength $ \gamma \in\R $.
To rigorously define the model,
we recall that $ x = (x_i)_{i=1}^\infty \in\R^{\N} $ is \textbf{rankable} 
if there exists a ranking permutation $ \per : \N \to \N $ 
such that $ x_{\per(i)} \leq x_{\per(j)} $, for all $ i < j \in \N $.
To ensure that such a ranking permutation is unique,
we resolve ties in lexicographic order.
That is, if $ x_{\per(i)} = x_{\per(j)} $ for $ i<j $, then $ \per(i) < \per(j) $.
We then let $ \per_x(\Cdot): \N \to \N $
denote the unique ranking permutation for a given, rankable $ x $.
Fix independent standard Brownian motions $ W_1,W_2,\ldots $. 
For suitable initial conditions,
the infinite Atlas model $ X(t) = (X_i(t))_{i=1}^\infty $
is given by the unique weak solution of the following system of \acp{SDE}
\begin{align}
	\label{eq:atlas}
	dX_i(t) = \gamma \ind\set{ \per_{X(t)}(i)=1 } dt + dW_i(t),
	\quad
	i\in\N.
\end{align}
To state the well-posedness results of \eqref{eq:atlas},
consider the following configuration space
\begin{align}
	\label{eq:cnfsp}
	\calU
	= 
	\Big\{ 
		x = (x_i)_{i=1}^\infty \, : \, 
		\lim_{i \to \infty}x_i = \infty, 
		\text{ and } \sum_{i=1}^{\infty}e^{-a x_i^2} < \infty,
		\forall a > 0
	\Big\},
\end{align}
and note that $ \lim_{i \to \infty}x_i = \infty $ necessarily implies that $ x $ is rankable.
It is shown in \cite[Theorem~3.2]{sarantsev14},
for any fixed $ \gamma\in\R $ and any given $ x\in\calU $,
the system~\eqref{eq:atlas} admits a unique weak solution $ X(t) $
starting from the initial condition $ x $,
such that $ \Pr(X(t)\in\calU,\ \forall t \geq 0 ) =1 $.
See also \cite{shkolnikov11,ichiba13}.

The interest of the Atlas model originates from the study 
of diffusions with rank-based drifts \cite{fernholz02,karatzas09}.
In particular, the Atlas model was first introduced, in finite dimensions,
as a simple special case of rank-based diffusions \cite{fernholz02}.
Due to their intriguing properties,
rank-based diffusions have been intensively studied in various generality. 
See \cite{banner05,banner11,ichiba13,sarantsev15} and the references therein.
The infinite-dimensional system~\eqref{eq:atlas} considered here
was introduced by Pal and Pitman~\cite{pal08}.
Parts of the motivation was to
understand the effect of a drift exerted on a large (but finite) collection of Brownian particles
\cite{aldous02,tang15}.
In particular, it was shown in \cite{pal08} that, for $ \gamma>0 $, 
the system~\eqref{eq:atlas}
admits a stationary gap distribution of i.i.d.\ $ \Exp(2\gamma) $,
which indicates that the drift $ \gamma dt $
is balanced by the push-back of a crowd of particles of density $ 2\gamma $.
To state the previous result more precisely,
given a rankable $ x=(x_i)_{i=1}^\infty $,
we let $ (x_{(1)}\leq x_{(2)}\leq\ldots) $ 
denote the corresponding ranked points, i.e., $ x_{(i)} = x_{(\per_{x})^{-1}(i)} $,
and consider the corresponding gaps $ z_i := x_{(i+1)}-x_{(i)} $.
It was shown in \cite{pal08} that $ \pi:=\bigotimes_{i=1}^\infty \Exp(2\gamma) $
is a stationary distribution of the gap process 
$ Z(t):= (X_{(i+1)}(t)-X_{(i)}(t))_{i=1}^\infty $ of the Atlas model~\eqref{eq:atlas}.

It addition to the i.i.d.\ $ \Exp(2\gamma) $ distribution,
it was recently shown in \cite{sarantsev17} that
the Atlas model has a different type of stationary gap distributions.
That is, for each $ a>2\gamma_- $, $ \pi_a := \bigotimes_{i=1}^\infty \Exp(2\gamma+ia) $
is also a stationary gap distribution of the Atlas model.
Unlike $ \pi $,
the distribution $ \pi_a $ has exponentially growing particle density away from the Atlas particle.
In this article, we go one step further and show that,
in fact, up to a deterministic shift $ \frac{at}{2} $ of each particle,
the \emph{entire} particle system 
$ \{ X_{i}(t)+\frac{at}{2} \}_{i=1}^\infty $
has a stationary distribution.
This extends the result of \cite{sarantsev17}
on stationary \emph{gap} distributions.
In the following we use $ \{x_{i}\}_{i=1}^\infty \subset\R $
to denote a configuration of \emph{indistinguishable} particles,
in contrast with $ (x_{i})_{i=1}^\infty $,
which denotes labeled (named) particles.
Let
\begin{align*}
	\calV
	= 
	\Big\{ 
		\{x_i\}_{i=1}^\infty \, : \, 
		(x_{i})_{i=1}^\infty \in \calU
	\Big\}
\end{align*}
denote the corresponding configuration space,
and let $ \poi_{a} $ denote the Poisson point process on $ \R $
with density $ ae^{a\xi}d\xi $.
It is standard to show (e.g., using techniques from~\cite[Section~2.2]{panchenko13})
that $ \poi_{a} $ is supported on $ \calV $.
Let $ \Gamma(\alpha) := \int_0^\infty \xi^{-1-\alpha} e^{-\xi} d\xi $
denote the Gamma function,
and let $ \Gammad(\alpha,\beta) \sim \frac{1}{\Gamma(\alpha)} \beta^\alpha\xi^{-1-\alpha} e^{-\beta\xi} \ind_\set{\xi >0} d\xi $
denote the Gamma distribution.
The following is the main result.

\begin{theorem}
\label{thm:main}
\begin{enumerate}[label=(\alph*)]
\item []
\item \label{enu:main1}
	For any fixed $ \gamma\in\R $ and  $ a >2\gamma_- $,
	$
		\Ex_{\poi_a}(e^{2\gamma X_{(1)}}) = \Gamma(\frac{2\gamma}{a}+1) \in (0,\infty),
	$
	so that
	\begin{align}
		\label{eq:invt}
		\invt_a(\Cdot) := \frac{1}{\Gamma(\frac{2\gamma}{a}+1)}\Ex_{\poi_{a}}(e^{2\gamma X_{(1)}}\Cdot)
	\end{align}
	defines a probability distribution supported on $ \calV $.
	Furthermore, under $ \invt_a $, we have that
	$ e^{aX_{(1)}}\sim\Gammad(\frac{2\gamma}{a}+1,1) $, and that
	\begin{align}
	\label{eq:renyi}
		Z := ( X_{(i+1)}-X_{(i)} )_{i=1}^\infty \sim \pi_a = \bigotimes_{i=1}^\infty \Exp(2\gamma+ia).
	\end{align}
\item	\label{enu:main2}
	The distribution $ \invt_a $ is a stationary distribution of
	$ \{ X_{i}(t)+\frac{at}{2} \}_{i=1}^\infty $,
	where $ (X_i(t))_{i=1}^\infty $ evolves under \eqref{eq:atlas}.
\end{enumerate}
\end{theorem}

\begin{remark}
Under $ \invt_a $, the Atlas particle $ X_{(1)} $ and the gap process $ Z=(Z_i)_{i=1}^\infty $ 
are \emph{not} independent.
\end{remark}

For the special case $ \gamma=0 $,
the Atlas model~\eqref{eq:atlas} reduces to independent Brownian motions.
In this case, it is well known that the Poisson point process $ \poi_{a} $ 
is quasi-stationary \cite{liggett78},
and the shift $ -\frac{a}{2}t $ can be easily calculated from 
the motion of independent Brownian particles. 
Here we show that, with a drift $ \gamma dt $ exerted on the Atlas particle $ X_{(1)}(t) $,
a stationary distribution is obtained by taking $ V(x) := 2\gamma x_{(1)} $ to be the potential.
Indeed, under such a choice of $ V(x) $, we have that
$ \gamma\ind\set{p_{x}(i)=1} = \tfrac12 \partial_{x_i} e^{V(x)} $.
This explains why we should expect the stationary distribution $ \invt_a $ as in \eqref{eq:invt}.
The proof of Theorem~\ref{thm:main} amounts to justifying 
the aforementioned heuristic in the setting of infinite dimensional diffusions
with discontinuous drift coefficients.
We achieve this through finite-dimensional, smooth approximations,
and using the explicit expressions of semigroups from Girsanov's theorem
to take limits.

Due to their simplicity,
product-of-exponential stationary gap distributions
have been intensively searched within competing Brownian particle systems,
in both finite and infinite dimensions.
See~\cite{sarantsev14} and the references therein.
To date, derivations of product-of-exponential stationary gap distributions
have been relying on the theory of Semimartingale Reflecting Brownian Motions (SRBM),
e.g., \cite{williams95}.
On the other hand,
given the expression~\eqref{eq:invt} of $ \invt_a $,
the gap distribution~\eqref{eq:renyi}
follows straightforwardly from R{\'e}nyi's representation~\cite{renyi53}.
Theorem~\ref{thm:main} hence provides an alternative 
derivation of the product-of-exponential distribution~$ \pi_a $
without going through SRBM.

Our methods should generalize to the case of competing Brownian particle systems
with finitely many non-zero drift coefficients, i.e.,
\begin{align*}
	dX_i(t) = \sum_{j=1}^m \gamma_j \ind\set{ \per_{X(t)}(i)=j } dt + dW_i(t),
	\quad
	i\in\N,
\end{align*}
yielding the stationary distribution 
$ \invt_a(\Cdot) := \frac{1}{J}\Ex_{\poi_{a}}(e^{2\sum_{j=1}^m \gamma_jX_{(j)}}\Cdot) $,
for some normalizing constant $ J<\infty $.
Here we consider only the Atlas model for simplicity of notations.

A natural question,
following the discovery a stationary gap distribution,
is the longtime behavior of the Atlas particle $ X_{(1)}(t) $
under such a gap distribution.
For the i.i.d.\ $ \Exp(2) $ gap distribution $ \pi $,
this question was raised in \cite{pal08} and answered in~\cite{dembo17}.
It was shown in~\cite{dembo17} that $ X_{(1)}(t) $ fluctuates at order $ t^{\frac14} $ 
around its starting location,
and scales to a $ \frac14 $-fractional Brownian motion, as $ t\to\infty $.
As a simple application of Theorem~\ref{thm:main},
under the stationary gap distribution $ \pi_a $
and $ X_{(1)}(0)=0 $,
we establish an exponential tail bound, uniformly in $ t $,
of the fluctuation Atlas particle around its expected location $ -\frac{at}{2} $.
This shows that the fluctuation of $ X_{(1)}(t) $ stays bounded
under $ \pi_a $,
in sharp contrast with the $ t^\frac14 $ fluctuation obtained in \cite{dembo17}.

\begin{corollary}
\label{cor}
Fix $ \gamma\in\R $ and $ a>2\gamma_- $.
Starting the Atlas model~\eqref{eq:atlas} from the initial distribution
$ X_{(i)}(0)=0 $, $ (X_{(i+1)}(0)-X_{(i)}(0))_{i=1}^\infty \sim \pi_a $,
we have that
\begin{align*}
	\Pr(|X_{(1)}(t)+\tfrac{at}{2}| \geq \xi)
	\leq
	c e^{-\frac12(2\gamma+a)\xi},
	\quad
	\forall t,\xi \in\R_+,
\end{align*}
for some constant $ c=c(a,\gamma)<\infty $
depending only on $ a,\gamma $.
\end{corollary}

\subsection*{Acknowledgements}
I thank Amir Dembo, Ioannis Karatzas, Andrey Sarantsev and Ofer Zeitouni
for enlightening discussions.
I am grateful to Andrey Sarantsev for many useful comments
on the presentation of this article.
This work was partially supported by a Junior Fellow award from the Simons Foundation to Li-Cheng Tsai.

\section{Proof}
\subsection{Theorem~\ref{thm:main}\ref{enu:main1}}
\label{sect:main1}
Fix $ \gamma\in\R $, $ a>2\gamma_- $,
and let $ \{X_i\}_{i=1}^\infty $ denote a sample from the Poisson point process $ \poi_{a} $.
Let $ N(\xi) := \#\{ X_i\in(-\infty,\xi] \} $ 
denote the number of particles in $ (-\infty,\xi] $,
whereby $ N(\xi) \sim \Pois(e^{a\xi}) $.
Indeed, $ \Pr_{\poi_{a}}(X_{(1)} >\xi) = \Pr(N(\xi)=0) = e^{-e^{a\xi}} $.
From this we calculate
\begin{align*}
	\Ex_{\poi_{a}}( e^{2\gamma X_{(1)}} \ind\set{ X_{(1)} \leq \xi } )
	=
	\int_{-\infty}^{\xi} e^{2\gamma \zeta} \frac{d~}{d\zeta} ( 1-  e^{-e^{a\zeta}}) d\zeta.
\end{align*}
Performing the change of variable $ \zeta' := e^{a\zeta} $, 
we see that\\
$ 
	\Ex_{\poi_{a}}( e^{2\gamma X_{(1)}} \ind\set{ e^{aX_{(1)}} \leq e^{\xi} } )
	=
	\int_{0}^{e^\xi} {\zeta'}^{\frac{2\gamma}{a}} e^{-\zeta'} d\zeta'.
$
From this it follows that 
$ \Ex_{\poi_{a}}( e^{2\gamma X_{(1)}} ) = \Gamma(\frac{2\gamma}{a}+1) $
and that $ e^{aX_{(1)}} \sim \Gammad(\frac{2\gamma}{a}+1,1) $
under $ \invt_a $.

Turning to showing~\eqref{eq:renyi},
we let $ \{X_i\}_{i=1}^\infty $ be sampled from $ \poi_{a} $
and let $ (Z_k)_{k=1}^\infty = (X_{(i+1)}-X_{(i)})_{i=1}^\infty $ denote the gap process.
Fix arbitrary $ m<\infty $.
Our goal is to show that
\begin{align}
	\label{eq:iidexp:goal}
	\frac{
		\Ex_{\poi_{a}}(e^{2\gamma X_{(1)}}\prod_{i=1}^{m-1} \ind\set{Z_i \geq \xi_i})
		}{
		\Ex_{\poi_{a}}(e^{2\gamma X_{(1)}})
	}
	=
	\prod_{i=1}^{m-1} e^{-(2\gamma+ia)\xi_i}
	=:
	\eta.
\end{align}
For any given threshold $ \xi\in\R $,
we let $ \poi_{a,\xi} $ denote the restriction of the Poisson point process $ \poi_a $ onto $ (-\infty,\xi] $.
For the restricted process,
we have $ \poi_{a,\xi} \sim \{ \xi-Y_1,\ldots, \xi-Y_{N(\zeta)} \} $,
where $ Y_1,Y_2,\ldots $ are i.i.d.\ $ \Exp(a) $ variables, independent of $ N(\xi) $.
Let $ Y_{(1)}< Y_{(2)} < \ldots < Y_{(n)} $ denote the ranking of $ (Y_1,\ldots,Y_n) $.
We then have that, conditionally on $ N(\xi) \geq m $,
$ (X_{(1)},\ldots,X_{(m)}) = (\xi-Y_{N(\xi)},\ldots,\xi-Y_{N(\xi)-m+1}) $.
Further, by R{\'e}nyi's representation~\cite{renyi53},
\begin{align*}
	(Y_{(k)})_{k=1}^n \stackrel{\text{d}}{=}
	\Big( \sum_{i=k}^{n}G_i \Big)_{k=1}^n,
	\quad
	\text{ where } \big(G_i \big)_{i=1}^n \sim \prod_{i=1}^n \Exp(ia).
\end{align*}
Using this we calculate
\begin{align*}
	\Ex_{\poi_{a}}\Big( e^{2\gamma X_{(1)}}&\prod_{i=1}^{m-1} \ind\set{Z_i \geq \xi_i} \Big| N(\xi)\geq m \Big)
\\
	&=
	\Ex\Big( 
		e^{2\gamma (N(\xi)\xi -\sum_{i=m}^{N(\xi)}G_i)} 
		\prod_{i=1}^{m-1} e^{-2\gamma G_i}\ind\set{G_i \geq \xi_i} 
	\Big| N(\xi)\geq m \Big)
\\
	&=
	\Ex\Big( e^{2\gamma (N(\xi)\xi -\sum_{i=m}^{N(\xi)}G_i}) \Big| N(\xi)\geq m \Big) \prod_{i=1}^{m-1} \frac{iae^{-(ai+2\gamma)\xi_i}}{ai+2\gamma}.
\end{align*}
Further use $ \frac{ia}{ai+2\gamma} = \Ex e^{-2\gamma G_i} $
to write $ \prod_{i=1}^{m-1} \frac{iae^{-(ai+2\gamma)\xi_i}}{ai+2\gamma} = \Ex(e^{-2\gamma(G_1+\ldots+G_{m-1})})\eta  $,
We then obtain
\begin{align*}
	\Ex_{\poi_{a}}\Big( e^{2\gamma X_{(1)}}\prod_{i=1}^{m-1} \ind\set{Z_i \geq \xi_i} \Big| N(\xi)\geq m \Big)
	=
	\Ex_{\poi_a}\Big( e^{2\gamma X_{(1)}} \Big| N(\xi)\geq m \Big) \eta.
\end{align*}
Taking into account the case $ N(\xi)<m $,
we write
\begin{align}
	\notag
	&\Ex_{\poi_{a}}\Big(e^{2\gamma X_{(1)}}\prod_{i=1}^{m-1} \ind\set{Z_i \geq \xi_i} \Big)
\\
	\label{eq:iidexp:}
	&=
	\Ex_{\poi_a}\Big( e^{2\gamma X_{(1)}} \ind\set{N(\xi)\geq m}\Big) \eta
	+
	\Ex_{\poi_a}\Big( e^{2\gamma X_{(1)}} \prod_{i=1}^{m-1} \ind\set{Z_i \geq \xi_i}) \ind\set{N(\xi)< m}\Big).
\end{align}
Since $ a>-2\gamma $,
fixing $ q>1 $ with $ |q-1| $ small enough,
we have 
\begin{align}
	\label{eq:unfInt}
	\Ex_{\poi_a}( e^{2q\gamma X_{(1)}}) = \Ex_{\poi_a}( |e^{2\gamma X_{(1)}}|^q) = \Gamma(\tfrac{2q\gamma}{a}+1)<\infty.
\end{align}
That is, $ e^{2q\gamma X_{(1)}} $ has bounded $ q $-th moment with $ q>1 $,
so in particular $ \{e^{2\gamma X_{(1)}}\ind\set{N(\xi)\geq m}\}_{\xi>0} $ is uniformly integrable.
For fixed $ m<\infty $, $ \ind\set{N(\xi)< m} \to_\text{P} 0 $, as $ \xi\to\infty $.
Using this to take the limit $ \xi\to\infty $ in \eqref{eq:iidexp:},
we thus obtain
$ 
	\Ex_{\poi_{a}}(e^{2\gamma X_{1}}\prod_{i=1}^{m-1} \ind\set{Z_i \geq \xi_i} ) 
	=  
	\Ex_{\poi_a}( e^{2\gamma X_{1}} ) \eta.
$
This concludes~\eqref{eq:iidexp:goal}.

\subsection{Theorem~\ref{thm:main}\ref{enu:main2}}
\label{sect:main2}
Samples from $ \poi_a $ have, almost surely,
no repeated points, i.e., $ X_{(1)}< X_{(2)} < X_{(3)}<\ldots $.
Fix arbitrary $ m<\infty $ and $ \phi\in C^\infty_c(\calW) $,
where $ \calW := \{ (x_1<x_2<\ldots<x_m) \} $ denote the Weyl chamber.
Let $ \bar{X}_i(t) := X_i(t) + \frac{a}{2}t $,
and $ \bar{X}_{(i)}(t) := X_{(i)}(t) + \frac{a}{2}t $ denote the compensated particle locations.
It then suffices to show that
\begin{align}
	\label{eq:invt:phi:}
	\Ex_{\poi_a}\big( e^{2\gamma \bar{X}_{(1)}(0)} \phi(\bar{X}_{(1)}(t),\ldots,\bar{X}_{(m)}(t)) \big)
	= \Ex_{\poi_a}\big( e^{2\gamma \bar{X}_{(1)}(0)} \phi(\bar{X}_{(1)}(0),\ldots,\bar{X}_{(m)}(0)) \big).
\end{align}
As will be convenient for notations,
for $ n\geq m $,
we consider the symmetric extension $ \phis $ of $ \phi $, defined 
for $ n\geq m $ as
\begin{align}
	\label{eq:phis}
	\phis: \R^n \to \R,
	\quad
	\phis(x) := 
		\phi(x_{(1)},\ldots,x_{(m)}).
\end{align}
We have slightly abused notations by using the same symbol $ \phis $ 
to denote the function for \emph{all} 
$ n\in \N_{\geq m}\cup\{\infty\} $.
Note that, by definition, the function $ \phi $ vanishes near the boundary
$ \{ (x_1 \leq \ldots x_i=x_{i+1} \ldots\leq x_m) : i=1,\ldots,n-1 \} $ of $ \calW $,
so, for $ n<\infty $, $ \phis\in \Cc^\infty(\R^n) $.

The strategy of proving~\eqref{eq:invt:phi:} is to approximate 
the infinite system~$ \bar{X}(t) $ by finite systems.
Fixing $ m\leq n<\infty $,
we consider the following $ n $-dimensional analog of~$ \bar{X}(t) $:
\begin{align}
	\label{eq:atlasn}
	\bar{X}^n_i(t) = x_i 
	+ \int_0^t \big( \gamma\ind\set{\per_{\bar{X}^n(s)}(i)=1} + \tfrac{a}{2} \big) ds 
	+ W_i(t),
	\quad
	i=1,\ldots,n,
\end{align}
where the ranking permutation $ \per_x(\Cdot) : \{1,\ldots,n\} \to \{1,\ldots,n\} $
is defined similarly to the case of infinite particles.
As the discontinuity of $ x\mapsto \ind\set{\per_{x}(i)=1} $
imposes unwanted complication in the subsequence analysis,
we consider further the mollified system as follows.
Fix a mollifier $ r\in C^\infty(\R^n) $,
i.e., $ r \geq 0 $, $ r|_{\Vert x\Vert\geq 1} =0 $ and $ \int_{\R^{n}} r(y)dy=1 $.
Let $ V(x) := 2\gamma x_{(1)} = 2\gamma \min(x_1,\ldots,x_n) $.
For $ \e \in(0,1) $, we define the mollified potential as
$ V^\e(x) := \int_{\R^n} V(y) r(\e^{-1}(x-y)) \e^{-n} dy $.
Under these notations, we have that
\begin{align}
	\label{eq:Venbh}
	\tfrac12\partial_i V^\e(x) = \gamma\ind\set{\per_{\bar{X}^n(s)}(i)=1},
	\quad \text{ on }
	\Omega_\e := \{ x\in\R^n : |x_i-x_j|>\e, \ \forall i<j  \}.
\end{align}
We then consider the following mollified system 
\begin{align}
	\label{eq:atlasne}
	\bar{X}^{n,\e}_i(t) = x_i 
	+ \int_0^t \big( \tfrac12\partial_{i}V^\e(\bar{X}^{n,\e}(s)) + \tfrac{a}{2} \big) ds 
	+ W_i(t),
	\quad
	i=1,\ldots,n.
\end{align}
With $ \partial_iV^\e $ being smooth and bounded,
the well-posedness of~\eqref{eq:atlasne} follows from standard theory, e.g., \cite{stroock07}.
Furthermore, letting $ u(t,x) := \Ex_x( \phis(\bar{X}^{n,\e}(t)) ) $,
we have that $ u\in C^{\infty}(\R_+\times\R^n) $,
and that $ u $ solves the following PDE:
\begin{align}
	\label{eq:uPDE}
	\partial_t u 
	= 
	\sum_{i=1}^n \big( \tfrac12 \partial_{ii} + \tfrac{a}{2} \partial_i + \tfrac12\partial_i V^\e \big) u,
	\quad
	u(0,x) = \phis(x).
\end{align}
With $ \partial_i V^\e $ being bounded and $ \phis $ being compactly supported,
applying the Feynman-Kac formula 
to the solution $ u $ of \eqref{eq:uPDE},
we see that $ u $ decays exponentially as $ |x|\to\infty $, i.e.,
\begin{align}
	\label{eq:utail:}
	\sup_{s\leq t, x\in\R^n} 
	\big\{ |u(t,x)| e^{\xi(|x_1|+\ldots+|x_n|)} \big\}
	<\infty,
	\quad
	\forall \xi,t<\infty.
\end{align}
Such an exponential 
estimate~\eqref{eq:utail:} progresses to higher order derivatives of $ u $.
More precisely,
with $ \partial_iV^\e \in C^\infty(\R^n) $ and $ u\in C^\infty(\R_+\times\R^n) $,
taking derivative $ \partial_{i} $ in \eqref{eq:uPDE},
we see that $ \partial_{i} u $ solves the following equation:
\begin{align}
	\label{eq:uPDE:}
	\partial_t (\partial_{i} u)
	=
	\sum_{j=1}^n
	\big(\big(
		\tfrac12 \partial_{jj}
		+ \tfrac{a}{2} \partial_j
		+ \partial_j V^\e
		\big)&(\partial_{i} u)
    	+(\partial_{ij}V^\e) u    
    \big),
\\
	&
    \notag
	\partial_{i} u(0,x) = \partial_{i} \phis(x) \in C^\infty_c(\R^n),
\end{align}
A similarly procedure applied to the solution $ \partial_i u $ of \eqref{eq:uPDE:}
yields
\begin{align*}
	\sup_{s\leq t, x\in\R^n, i=1,\ldots,n} 
	\big\{ |\partial_{i}u(t,x)| e^{\xi(|x_1|+\ldots+|x_n|)} \big\}
	<\infty,
	\quad
	\forall \xi,t<\infty.
\end{align*}
Iterating this argument to higher order derivatives, we obtain
\begin{align}
	\label{eq:utail}
	\sup_{s\leq t, x\in\R^n,|\beta|\leq k} 
	\big\{ |\partial_{\beta} u(t,x)| e^{\xi(|x_1|+\ldots+|x_n|)} \big\}
	<\infty,
	\quad
	\forall \xi,t,k<\infty.
\end{align}
The PDE~\eqref{eq:uPDE} has stationary distribution
$ e^{V^\e(x)} \prod_{i=1}^n e^{ax_i}dx_i $
(\emph{not} a probability distribution, since the total mass is infinite).
More precisely, integrate $ u(t,x) $ against the aforementioned distribution
to get 
\begin{align*}
	v(t) := \int_{\R^n} u(t,x)e^{V^\e(x)} \prod_{i=1}^n e^{ax_i}dx_i.
\end{align*}
Taking time derivative using \eqref{eq:uPDE} and \eqref{eq:utail},
followed by integrations by parts
\begin{align*}
	\int_{\R^n} \tfrac12(\partial_{ii} u(t,x))e^{V^\e(x)} \prod_{j=1}^n e^{ax_j}dx_j
	=
	-\int_{\R^n} (\partial_{i} u(t,x)) (\tfrac12\partial_iV^\e(x)+\tfrac{a}{2}) 
	e^{V^\e(x)} \prod_{j=1}^n e^{ax_j}dx_j,
\end{align*}
we obtain that $ \frac{d~}{dt} v(t) =0 $.
Consequently,
\begin{align}
	\label{eq:invt:ne}
	\int_{\R^n}  \Ex_{x}\big( \phis(\bar{X}^{n,\e}(t)) \big) e^{V^\e(x)} \prod_{i=1}^n e^{ax_i}dx_i
	=
	\int_{\R^n}  \phis(x) e^{V^\e(x)} \prod_{i=1}^n e^{ax_i}dx_i.
\end{align}

The next step is to take the limit $ \e\to 0 $ in \eqref{eq:invt:ne},
for \emph{fixed} $ n $.
This amounts to establishing the convergence of the term $ \Ex_{x}\big( \phis(\bar{X}^{n,\e}(t)) \big) $.
To this end, we use Girsanov's theorem to write
\begin{align}
	\label{eq:grsnv:n}
	\Ex_{x}\big( \phis(\bar{X}^{n}(t)) \big)
	&=
	\Ex_x \big( \phis(H(t)) F(t) \big),
\\
	\label{eq:grsnv:ne}
	\Ex_{x}\big( \phis(\bar{X}^{n,\e}(t)) \big)
	&=
	\Ex_x \big( \phis(H(t)) F^\e(t) \big),
\end{align}
where $ H(t) := (W_i(t)+\frac{at}{2}+x_i)_{i=1}^n $
consists of independent, drifted Brownian motions starting from $ x=(x_i)_{i=1}^n $,
and the terms $ F(t) $ and $ F^\e(t) $ are stochastic exponentials given by
\begin{align}
	\label{eq:F}
	F(t) &:=
	\exp\big( M(t) - \tfrac12 \langle M\rangle(t) \big),&
	&
	M(t) := \int_0^t\sum_{i=1}^n \gamma \ind\set{ \per_{H(t)}(i)=1 }dW_i(s),
\\
	\label{eq:Fe}
	F^\e(t) &:=
	\exp\big( M^\e(t) - \tfrac12 \langle M^\e\rangle(t) \big),&
	&
	M^\e(t) := 
		\int_0^t\sum_{i=1}^n \frac12\partial_i V^\e(H(s)) dW_i(s).
\end{align}
Taking the difference of \eqref{eq:grsnv:n}--\eqref{eq:grsnv:ne},
followed by using the Cauchy--Schwarz inequality, we obtain
\begin{align}
	\notag
	\big| 
		&\Ex_{x}\big( \phis(\bar{X}^{n}(t)) \big) 
		-\Ex_{x}\big( \phis(\bar{X}^{n,\e}(t)) \big)
	\big|
	=
	\big| \Ex_x \big( \phis(H(t)) F(t) (1-\tfrac{F^\e(t)}{F(t)}) \big) |
\\
	\label{eq:cauchy}
	&\leq 
	\Big( \Ex_x \big( \phis(H(t))^2 F(t)^2 \big) \Big)^\frac12
	\Big( \Ex_x (1-\tfrac{F^\e(t)}{F(t)})^2 \Big)^\frac12.
\end{align}
For the two terms in \eqref{eq:cauchy},
we next show that: \textit{i}) the first term is bounded; and 
\textit{ii}) the second term vanishes as $ \e\to 0 $.
Hereafter, we use $ c(a_1,a_2,\ldots) $
to denote a finite, deterministic constant,
that may change from line to line,
but depends only on the designated variables $ a_1,a_2,\ldots $.

\textit{i})
Recall that $ \phis $ is defined in terms of $ \phi $ through~\eqref{eq:phis}.
We fix $ \lambda<\infty $, \emph{independently} of $ n $, 
such that $ \supp(\phis) \subset [-\lambda,\lambda]^n $.
Under these notations,
\begin{align}
	\notag
	\Ex_x \big( \phis(H(t))^2 F(t)^2 \big) 
	&\leq
	\Vert \phi \Vert^2_{L^\infty} 
	\Ex_x \big( \ind_\set{H(t)\in[-\lambda,\lambda]^n} F(t)^2 \big) 
\\
	\label{eq:bd1}
	&\leq
	\Vert \phi \Vert^2_{L^\infty} 
	\big( \Ex F(t)^4 \big)^\frac12 \Pr_x \big( H(t)\in[-\lambda,\lambda]^n \big)^\frac12.
\end{align}
With $ F(t) $ defined in \eqref{eq:F},
and $ \langle M \rangle(t) = \gamma^2t $,
it follows that
\begin{align*}
	\Ex_x(F(t)^4) 
	=
	\Ex_x \big(e^{4M(t)} e^{ -2\langle M\rangle(t)} \big)
	= e^{\frac12(16-4)\gamma^2t} = c(\gamma,t).
\end{align*}
Let $ \Phi(x) := \int_{-\infty}^{x} \frac{1}{\sqrt{2\pi}}e^{-\frac{y^2}{2}} dy $ denote the Gaussian distribution function.
With $ H_i(t)= x_i+\frac{a}{2}t + W_i(t) $, we have
\begin{align*}
	\Pr_x \big( H(t)\in[-\lambda,\lambda]^n \big) 
	\leq 
	\prod_{i=1}^n \Phi\Big( \frac{\lambda -\tfrac{a}{2}t - x_i}{\sqrt{t}} \Big).
\end{align*}
Inserting these bounds into~\eqref{eq:bd1},
we obtain
\begin{align}
	\label{eq:1st:bd}
	\Ex_x \big( \phis(H(t))^2 F(t)^2 \big) 
	&\leq
	c(a,\gamma,\lambda,t) \prod_{i=1}^n \Phi\Big( \frac{\lambda -\tfrac{a}{2}t - x_i}{\sqrt{t}} \Big)
\\
	\label{eq:1st:bd:}
	&\leq
		c(a,\gamma,\lambda,t,n) \exp\Big(-\frac{x_1^2+\ldots+x_n^2}{4(t+1)}\Big).
\end{align}

\textit{ii})
Expand the expression $ \Ex_x (1-\tfrac{F^\e(t)}{F(t)})^2 $ into
\begin{align}
	\label{eq:FFexpnd}
	\Ex_x (1-\tfrac{F^\e(t)}{F(t)})^2
	=
	1 + \Ex_x (\tfrac{F^\e(t)}{F(t)})^2 - 2 \Ex_x\tfrac{F^\e(t)}{F(t)}.
\end{align}
From \eqref{eq:F}--\eqref{eq:Fe},
we have
\begin{align}
	\label{eq:FfF}
	\tfrac{F^\e(t)}{F(t)} = \exp( M(t)-M^\e(t) ) \exp( \tfrac12\langle M\rangle(t)-\tfrac12\langle M^\e \rangle(t)).
\end{align}
Set $ b^\e_i(s) := \frac12\partial_i V^\e(H(s)) $ to simplify notations.
As $ V(x) $ is Lipschitz with Lipschitz seminorm $ 2|\gamma| $, 
(i.e., $ |V(x) - V(y)| \leq 2\gamma |x-y| $, $ \forall x,y \in \R^n $),
we have $ |b^\e_i(s)| \leq |\gamma| $.
Consequently,
\begin{align}
	\label{eq:qv:bd}
	\langle M \rangle(t) = \gamma^2 t,
	\quad
	\langle M^\e \rangle(t) \leq n\gamma^2 t.
\end{align}
To estimate the expression~\eqref{eq:FfF},
we use \eqref{eq:Venbh} and $ |b^\e_i(s)| \leq |\gamma| $ to write 
\begin{align}
	\label{eq:qd:bd1}
	|\langle M-M^\e \rangle(t)| 
	=
	\int_0^t \sum_{i=1}^n \big( b^\e_i(s) - \gamma\ind\set{ \per_{H(s)}(i)=1 } \big)^2 ds 
	\leq 
	4n\gamma^2
	\int_0^t \ind\set{ H(s) \notin \Omega_\e} ds .
\end{align}
Let $ L_{i,j}(s,\xi) $ denote the localtime process of $ H_i(s)-H_j(s) = W_i(s)-W_j(s)+(x_j-x_i) $ 
at a given level $ \xi $.
We further bound the r.h.s.\ of \eqref{eq:qd:bd1} as
\begin{align*}
	\int_0^t \ind\set{ H(s) \notin \Omega_\e} ds
	\leq
	\sum_{i<j} \int_0^t \int_{|\xi|\leq \e} L_{i,j}(s,\xi) d\xi ds
	\longrightarrow_\text{P}
	0,
	\quad
	\text{ as } \e \to 0.
\end{align*}
Consequently, $ |\langle M-M^\e \rangle(t)|\to_\text{P} 0 $.
Since, by \eqref{eq:qv:bd}, 
$ \langle M\rangle (t) $ and $ \langle M^\e\rangle (t) $ are bounded
(for fixed $ t $),
it also follows that $ \Ex_x|\langle M-M^\e \rangle(t)| \to 0 $
and hence $  M(t)-M^\e(t) \to_\text{P} 0 $.
Referring back to the expression~\eqref{eq:FfF},
we see that $ \frac{F^\e(t)}{F(t)}\to_\text{P} 1 $.
Using again the fact that 
$ \langle M\rangle (t) $ and $ \langle M^\e\rangle (t) $ are bounded,
(which implies the uniform integrability of $ (\frac{F^\e(t)}{F(t)})^k $, $ k=1,2$),
we obtain $ \Ex_x(\frac{F^\e(t)}{F}), \Ex_x(\frac{F^\e(t)}{F})^2\to 1 $.
Inserting these into \eqref{eq:FFexpnd} yields
\begin{align}
	\label{eq:2nd:bd}
	\Ex_x (1-\tfrac{F^\e(t)}{F(t)})^2 \to 0,
	\text{ as } \e \to 0,
	\quad
	\text{ for any fixed } x\in\R^n.
\end{align}

Now, combine \eqref{eq:1st:bd:}, \eqref{eq:2nd:bd} with \eqref{eq:cauchy},
and insert the result into the l.h.s.\ of \eqref{eq:invt:ne}.
After taking the $ \e\to 0 $ limit with $ n<\infty $ being fixed,
we obtain
\begin{align}
	\label{eq:invt:n}
	\int_{\R^n}  \Ex_{x}\big( \phis(\bar{X}^{n}(t)) \big) e^{2\gamma x_{(1)}} \prod_{i=1}^n e^{ax_i}dx_i
	=
	\int_{\R^n}  \phis(x) e^{2\gamma x_{(1)}} \prod_{i=1}^n e^{ax_i}dx_i.
\end{align}

Recall that $ \poi_{a,\zeta} $ denote the restriction of
the Poisson point process $ \poi_a $ on $ (-\infty,\zeta] $
and that $ N(\zeta) $ denote the number of particles on $ (-\infty,\zeta] $.
As mentioned previously,
$ \poi_{a,\zeta} \sim \{ \zeta-Y_1,\ldots, \zeta-Y_{N(\zeta)} \} $,
where
$ Y_1,Y_2,\ldots $ are i.i.d.\ $ \Exp(a) $ variables, independent of $ N(\zeta) $.
Conditionally on $ N(\zeta)=n $,
the process $ \{ \zeta-Y_1,\ldots, \zeta-Y_{N(\zeta)} \} $
have joint distribution $ \prod_{i=1}^n ae^{a(x_i-\zeta)}dx_i \ind_\set{x_i\leq\zeta} $.
With this, multiplying both sides of \eqref{eq:invt:n} by $ a^ne^{-an\zeta} $,
and averaging over $ \{N(\zeta) \geq m\} $,
we obtain that
\begin{align}
	\notag
	\Ex_{\poi_{a,\zeta}}&\Big( \phis(\bar{X}^{N(\zeta)}(t)) e^{2\gamma \bar{X}^{N(\zeta)}_{(1)}(0)} 
		\ind\set{N(\zeta) \geq m} 
	\Big)
	+ \Ex(R_{N(\zeta)}(\zeta)\ind\set{N(\zeta) \geq m} )	
\\
	\label{eq:invt:n:}
	&=
	\Ex_{\poi_{a,\zeta}}\Big( \phis(\bar{X}^{N(\zeta)}(0)) e^{2\gamma \bar{X}^{N(\zeta)}_{(1)}(0)} 
	\ind\set{N(\zeta) \geq m} \Big)
	+ \Ex(S_{N(\zeta)}(\zeta)\ind\set{N(\zeta) \geq m} ),
\end{align}

where the terms $ R_n(\zeta) $ and $ S_n(\zeta) $ are given by
\begin{align}
	\label{eq:R}
	R_n(\zeta) &:= \int_{\cup_{i=1}^n\{x_i >\zeta\}}  
	\Ex_{x}\big( \phis(\bar{X}^{n}(t)) \big) e^{2\gamma x_{(1)}} \prod_{i=1}^n a e^{a(x_i-\zeta)}dx_i,
\\
	\notag
	S_n(\zeta) &:= \int_{\cup_{i=1}^n\{x_i >\zeta\}}  
	\Ex_{x}\big( \phis(x) \big) e^{2\gamma x_{(1)}} \prod_{i=1}^n a e^{a(x_i-\zeta)}dx_i.
\end{align}
Recall that $ \supp(\phis)\subset[-\lambda,\lambda]^n $.
Hence
\begin{align}
	\label{eq:S:bd}
	S_{n}(\zeta) =0, \quad\text{for all } \zeta > \lambda.
\end{align}
As for $ R_n(\zeta) $,
inserting the bound~\eqref{eq:1st:bd} into \eqref{eq:R} gives
\begin{align*}
	|R_n(\zeta)|
	\leq
	c(a,\gamma,\lambda,t) 
	\int_{\cup_{i=1}^n\{x_i >\zeta\}}  
	e^{2\gamma x_{(1)}} \prod_{i=1}^n \Phi\Big( \frac{\lambda -\tfrac{a}{2}t - x_i}{\sqrt{t}} \Big) ae^{a(x_i-\zeta)}dx_i.
\end{align*}
Indeed, $ x_{(1)} \leq \zeta + \sum_{i=1}^n (x_i-\zeta)_+ $, so, after a change of variable $ x_i-\zeta\mapsto x_i $, we obtain
\begin{align*}
	|R_n(\zeta)|
	\leq
	c(a,\gamma,\lambda,t) e^{\zeta}
	\int_{\cup_{i=1}^n\{x_i >0\}}  
	\prod_{i=1}^n \Phi\Big( \frac{\lambda -\tfrac{a}{2}t - x_i-\zeta}{\sqrt{t}} \Big) ae^{ax_i+a (x_i)_+}dx_i.
\end{align*}
To bound the last integral, we split the integration over $ x_i $ into $ \{x_i > 0\} $ and $ \{x_i\leq 0\} $ for \emph{each} $ x_i $,
and thereby express the integral as
\begin{align*}
	\sum_{k=1}^n \sum_{\{i_1,\ldots,i_k\}} 
	\Big( 
		\prod_{j\in\{i_1,\ldots,i_k\}} \int_{\{x_j>0\}}(\cdots) dx_j
	\Big)
	\Big( 
		\prod_{j\notin\{i_1,\ldots,i_k\}} \int_{\{x_j\leq 0\}}(\cdots) dx_j
	\Big),
\end{align*}
where $ \{i_1,\ldots,i_k\} $ ranges over all distinct $ k $-indices from $ \{1,\ldots,n\} $.
Further, for each integral over $ \{x>0\} $ and over $ \{x \leq 0\} $, we have that
\begin{align*}
	&\int_{\{x > 0\}} \Phi\Big( \frac{\lambda -\tfrac{a}{2}t - x-\zeta}{\sqrt{t}} \Big) ae^{ax+a (x)_+}dx
	\leq
	c(a,\lambda,\gamma,t) e^{-\frac{\zeta^2}{4(t+1)}},
\\
	&\int_{\{x \leq 0\}} \Phi\Big( \frac{\lambda -\tfrac{a}{2}t - x-\zeta}{\sqrt{t}} \Big) ae^{ax+a (x)_+}dx
	<
	\int_{\{x \leq 0\}} a e^{ax} dx =1.
\end{align*}
Consequently,
\begin{align*}
	|R_n(\zeta)|
	\leq
	c(a,\gamma,\lambda,t) e^{\zeta}	\sum_{k=1}^n \binom{n}{k} c(a,\gamma,\lambda,t)^{k} e^{-\frac{k\zeta^2}{4(t+1)}}.  
\end{align*}
Now, with $ N(\zeta)\sim \Pois(e^{a\zeta}) $, we have $ \Ex(\binom{N(\zeta)}{k})= \frac{1}{k!}\Ex(N(\zeta)\cdots(N(\zeta)-k+1)) = \frac{1}{k!}e^{ka\zeta} $.
Given this identity, setting $ n=N(\zeta) $ and taking expected value, we obtain
\begin{align}
	\label{eq:R:bd}
	\Ex|R_{N(\zeta)}(\zeta)|
	\leq
	c(a,\gamma,\lambda,t) e^{\zeta}	\sum_{k=1}^\infty \frac{1}{k!} c(a,\lambda,t)^{k} e^{ka\zeta-\frac{k\zeta^2}{4(t+1)}},
\end{align}
which converges to zero as $ \zeta \to \infty $.

Using \eqref{eq:S:bd}--\eqref{eq:R:bd} in \eqref{eq:invt:n:},
and taking the limit $ \zeta\to \infty $, we arrive at
\begin{align}
	\notag
	\lim_{\zeta\to\infty}
	&\Big(
		\Ex_{\poi_{a,\zeta}}\Big( \phis(\bar{X}^{N(\zeta)}(t)) e^{2\gamma \bar{X}^{N(\zeta)}_{(1)}(0)} \ind\set{N(\zeta) \geq m}\Big)
\\
	\label{eq:invt:n::}
		&-
		\Ex_{\poi_{a,\zeta}}\Big( \phis(\bar{X}^{N(\zeta)}(0)) e^{2\gamma \bar{X}^{N(\zeta)}_{(1)}(0)} \ind\set{N(\zeta) \geq m}\Big)
	\Big)
	=0.
\end{align}
It remains to show that,
under the limit $ \zeta\to\infty $,
we can exchange the finite system~$ \bar{X}^{N(\zeta)} $ for the infinite system~$ \bar{X} $
within the expressions in \eqref{eq:invt:n::}.
As $ \zeta\to\infty $, we have that
\begin{align}
	\label{eq:XnX}
	\bar{X}^{N(\zeta)}_{(i)}(t) \Rightarrow \bar{X}_{(i)}(t),
	\
	\text{ as } \zeta\to\infty,
	\quad
	i=1,\ldots,m,
\end{align}
where $ \bar{X}^{N(\zeta)}(0) \sim \poi_{a,\zeta} $
and $ \bar{X}(0) \sim \poi_a $.
Such a statement~\eqref{eq:XnX} can be proven by
techniques from \cite{sarantsev14} and \cite[Section~3(a)]{sarantsev17}.
We omit repeating the standard arguments here.
Combining \eqref{eq:XnX} and \eqref{eq:unfInt}, we obtain that
\begin{align}
	\label{eq:invt:n::1}
	&
	\lim_{\zeta\to\infty}
	\Ex_{\poi_{a,\zeta}}\Big( \phis(\bar{X}^{N(\zeta)}(t)) e^{2\gamma \bar{X}^{N(\zeta)}_{(1)}(0)} \Big)
	=
	\Ex_{\poi_{a}}\Big( \phis(\bar{X}(t)) e^{2\gamma \bar{X}_{(1)}(0)} \Big),
\\
	\label{eq:invt:n::2}
	&
	\lim_{\zeta\to\infty}
	\Ex_{\poi_{a,\zeta}}\Big( \phis(\bar{X}^{N(\zeta)}(0)) e^{2\gamma \bar{X}^{N(\zeta)}_{(1)}(0)} \Big)
	=
	\Ex_{\poi_{a}}\Big( \phis(\bar{X}(0)) e^{2\gamma \bar{X}_{(1)}(0)} \Big).
\end{align}
Combining~\eqref{eq:invt:n::1}--\eqref{eq:invt:n::2} with \eqref{eq:invt:n::},
we thus obtain~\eqref{eq:invt:phi:},
and hence complete the proof.
\subsection{Corollary~\ref{cor}}
Fixing $ \gamma\in\R $ and $ a> 2\gamma_- $,
we let $ c=c(a,\gamma)<\infty $
denote a generic finite constant that depends only on these two variables.
Let $ Y(t) = (Y_i(t))_{i=1}^\infty $ be a solution to~\eqref{eq:atlas}
starting from the distribution $ \{Y_i(0)\}_{i=1}^\infty \sim \invt_a $,
so that $ \{Y_i(t)+\frac{at}{2}\}_{i=1}^\infty \sim \invt_a $,
for all $ t\in\R_+ $.
Since, by \eqref{eq:renyi}, the gap process $ (Y_{(i+1)}(0)-Y_{(i)}(0))_{i=1}^\infty $
is distributed as $ \pi_a $,
setting $ X_i(t)=Y_i(t)-Y_{(1)}(0) $,
we have that $ X(t) $ is a solution to~\eqref{eq:atlas} 
with the designated initial distribution as in Corollary~\ref{cor}.
Under these notations, for any given $ \xi \geq 0 $,
\begin{align}
	\notag
	\Pr( |X_{(1)}(t)| \geq \xi )
	&=
	\Pr( |Y_{(1)}(t)-Y_{(1)}(0)| \geq \xi )
\\
	\label{eq:Y1bd}
	&\leq
	\Pr( |Y_{(1)}(0)| \geq \tfrac{\xi}{2} ) + \Pr( |Y_{(1)}(t)| \geq \tfrac{\xi}{2} )
	=
	2\Pr( |Y_{(1)}(0)| \geq \tfrac{\xi}{2} ).
\end{align}
With $ e^{a Y_{(1)}} \sim \Gammad(\frac{2\gamma}{a},1) $,
we have that
\begin{align*}
	\Pr( Y_{(1)}(0) \leq -\tfrac{\xi}{2} )
	&=
	\frac{1}{\Gamma(\frac{2\gamma}{a})}\int_{0}^{e^{-\frac12a\xi}} \zeta^{\frac{2\gamma}{a}} e^{-\zeta} d\zeta
	\leq
	c\int_{0}^{e^{-\frac12a\xi}} \zeta^{\frac{2\gamma}{a}} d\zeta
	=
	c e^{-\frac12(2\gamma+a)\xi},
\\
	\Pr( Y_{(1)}(0) \geq \tfrac{\xi}{2} )
	&=
	\frac{1}{\Gamma(\frac{2\gamma}{a})}\int_{e^{\frac12a\xi}}^\infty \zeta^{\frac{2\gamma}{a}} e^{-\zeta} d\zeta
	\leq
	c\int_{e^{\frac12a\xi}}^\infty  e^{-\frac12\zeta} d\zeta
	=
	c e^{-\frac12e^{\frac12a\xi}}
	\leq
	c e^{-\frac12(2\gamma+a)\xi}.
\end{align*}
Combining these bounds with \eqref{eq:Y1bd}
yields the desired result.

\bibliographystyle{alphaabbr}
\bibliography{2na}

\end{document}